\documentclass[MSNbibl,number,citesort,dvips]{arxbj}

\aid{0}
\volume{20}
\issue{1}
\pubyear{2014}
\firstpage{282}
\lastpage{303}
\doi{10.3150/12-BEJ486} 

\makeatletter

\newcommand{\rrVert}{\Vert}
\newcommand{\rrvert}{\vert}
\newcommand{\llVert}{\Vert}
\newcommand{\llvert}{\vert}

\newtheorem{theorem}{Theorem}
\newtheorem{lemma}[theorem]{Lemma}
\newtheorem{proposition}[theorem]{Proposition}

\newproclaim{assumption}[theorem]{Assumption}

\newremark{Remark}{Remark}
\newremark{Remarks}{Remarks}

\def\rank{\operatorname{rank}}
\def\bR{\mathbb R}
\def\bE{\mathbb E}

\makeatother

\begin{document}
\begin{frontmatter}

\title{Noisy low-rank matrix completion with general sampling distribution}
\runtitle{Noisy low-rank matrix completion}

\begin{aug}
\author{\fnms{Olga} \snm{Klopp}\corref{}\ead[label=e1]{kloppolga@math.cnrs.fr}}
\runauthor{O. Klopp} 
\address{MODAL'X, University Paris Ouest Nanterre and CREST,
200 avenue de la R{\'e}publique, 92001 Nanterre, France. \printead{e1}}
\end{aug}

\received{\smonth{3} \syear{2012}}
\revised{\smonth{7} \syear{2012}}

%
\begin{abstract}
In the present paper, we consider the problem of matrix completion with
noise. Unlike previous works, we consider quite general sampling
distribution and we do not need to know or to estimate the variance of
the noise. Two new nuclear-norm penalized estimators are proposed, one
of them of ``square-root'' type. We analyse their performance under
high-dimensional scaling and provide non-asymptotic bounds on the
Frobenius norm error. Up to a logarithmic factor, these performance
guarantees are minimax optimal in a number of circumstances.
\end{abstract}

%
\begin{keyword}
\kwd{high-dimensional sparse model}
\kwd{low rank matrix estimation}
\kwd{matrix completion}
\kwd{unknown variance}
\end{keyword}

\end{frontmatter}

\section{Introduction}\label{introduction}
This paper considers the problem of matrix recovery from a small set of
noisy observations. Suppose that we observe a small set of entries of a
matrix. The problem of inferring the many missing entries from this set
of observations is the \textit{matrix completion} problem. A usual
assumption that allows to succeed such a completion is to suppose that
the unknown matrix has low rank or has approximately low rank.

The problem of matrix completion comes up in many areas including
collaborative filtering, multi-class learning in data analysis, system
identification in control, global positioning from partial distance
information and computer vision, to mention some of them. For instance,
in computer vision, this problem arises as many pixels may be missing
in digital images. In collaborative filtering, one wants to make
automatic predictions about the preferences of a user by collecting
information from many users. So, we have a data matrix where rows are
users and columns are items. For each user, we have a partial list of
his preferences. We would like to predict the missing rates in order to
be able to recommend items that may interest each user.

The noiseless setting was first studied by
Cand{\`e}s and Recht \cite{candes-recht-exact} using nuclear norm minimization.
A tighter analysis of the same convex relaxation was carried out in
\cite{candes-tao-power}. For a simpler approach, see more recent
papers of Recht \cite{recht-simpler} and Gross \cite{gross-recovery}.
An alternative line of work was developed by Keshavan \textit{et al.}  in \cite
{keshavan-few}.
A more common situation in applications corresponds to the noisy
setting in which the few available entries are corrupted by noise. This
problem has been extensively studied recently. The most popular methods
rely on nuclear norm minimization (see, e.g., \cite
{candes-plan-noise,keshavan-montanari-matrix,rhode-tsybakov-estimation,wainwright-weighted,koltchinskii-von,Koltchinskii-Tsybakov,2010arxiv10084886G,klopp-variance,foygel-serebro-concentration}).
One can also use rank penalization as it was done by Bunea \textit{et al.}  \cite
{bunea} and Klopp \cite{klopp-rank}.
Typically, in the matrix completion problem, the sampling scheme is
supposed to be uniform.
However, in practice, the observed entries are not guaranteed to follow
the uniform scheme and its distribution is not known exactly.

In the present paper, we consider nuclear norm penalized estimators and
study the corresponding estimation error in Frobenius norm. We consider
both cases when the variance of the noise is known or not. Our methods
allow us to consider quite general sampling distribution: we only
assume that the sampling distribution satisfies some mild
``regularity'' conditions (see Assumptions~\ref{L} and \ref{assPi}).

Let $A_0\in\mathbb R^{m_1\times m_2}$ be the unknown matrix. Our main
results, Theorems \ref{thmu3} and \ref{thm3}, show the following
bound on the normalized Frobenius error of the estimators $\hat A$ that
we propose in this paper: with high probability
\[
\frac{\Vert\hat A-A_0\Vert_2^{2}}{m_1m_2}\lesssim\frac{\log
(m_1+m_2)\max(m_1,m_2)\rank(A_0)}{n},
\]
where the symbol $\lesssim$ means that the inequality holds up to a
multiplicative numerical constant. This theorem guarantees, that the
prediction error of our estimator is small whenever $n\gtrsim\log
(m_1+m_2)\max(m_1,m_2)\rank(A_0)$. This quantifies the sample size
necessary for successful matrix completion.
Note that, when $\rank(A_0)$ is small, this is considerably smaller
than $m_1m_2$, the total number of entries. For large $m_1,m_2$ and
small $r$, this is also quite close to the degree of freedom of a rank
$r$ matrix, which is $(m_1+m_2)r-r^{2}$.

An important feature of our estimator is that its construction requires
only an upper bound on the maximum absolute value of the entries of
$A_0$. This condition is very mild. A~bound on the maximum of the
elements is often known in applications.
For instance, if the entries of $A_0$ are some user's ratings it
corresponds to the maximal rating. Previously, the estimators proposed
by Koltchinskii \textit{et al.}  \cite{Koltchinskii-Tsybakov} and by Klopp \cite
{klopp-rank} also require a bound on the maximum of the elements of the
unknown matrix but their constructions use the uniform sampling and
additionally require the knowledge of an upper bound on the variance of
the noise. Other works on matrix completion require more involved
conditions on the unknown matrix. For more details, see Section \ref
{completionknown}.

Sampling schemes more general than the uniform one were previously
considered in \mbox{\cite
{lounici-spectral,wainwright-weighted,serbro-learning}}. Lounici \cite
{lounici-spectral} considers a different estimator and measures the
prediction error in the spectral norm. In \cite
{wainwright-weighted,serbro-learning} the authors consider penalization
using a weighted trace-norm, which was first introduced by Srebro
\textit{et al.}  \cite{serbro-collaborative}. Negahban \textit{et al.}
in \cite {wainwright-weighted} assume that the sampling distribution is
a product distribution, that is, the row index and the column index of
the observed entries are selected independently. This assumption does
not seem realistic in many cases (see discussion in~\cite
{serbro-learning}). An important advantage of our method is that the
sampling distribution does not need to be equal to a product
distribution. Foygel \textit{et al.}  in \cite{serbro-learning} propose
a method based on the ``smoothing'' of the sampling distribution. This
procedure may be applied to an arbitrary sampling distribution but
requires a priori information on the rank of the unknown matrix.
Moreover, unlike in the present paper, in \cite{serbro-learning} the
prediction performances of the estimator are evaluated through a bound
on the expected $l$-Lipschitz loss (where the expectation is taken with
respect to the sampling distribution).

The weighted trace-norm, used in \cite
{wainwright-weighted,serbro-learning}, corrects a specific situation
where the standard trace-norm fails. This situation corresponds to a
non-uniform distribution where the row/column marginal distribution is
such that some columns or rows are sampled with very high probability
(for a more thorough discussion see \cite
{serbro-collaborative,serbro-learning}). Unlike \cite
{wainwright-weighted,serbro-learning}, we use the standard trace-norm
penalization and our assumption on the sampling distribution
(Assumption \ref{L}) guarantees that no row or column is sampled with
very high probability.

Most of the existing methods of matrix completion rely on the knowledge
or a pre-estimation of the standard deviation of the noise. The matrix
completion problem with unknown variance of the noise was previously
considered in \cite{klopp-variance} using a different estimator which
requires uniform sampling. Note also that in \cite{klopp-variance} the
bound on the prediction error is obtained under some additional
condition on the rank and the ``spikiness ratio'' of the matrix. The
construction of the present paper is valid for more general sampling
distributions and does not require such an extra condition.

The remainder of this paper is organized as follows. In Section \ref
{Preliminaries}, we introduce our model and the assumptions on the
sampling scheme. For the reader's convenience, we also collect notation
which we use throughout the paper. In Section \ref{completionknown}
we consider matrix completion in the case of known variance of the
noise. We define our estimator and prove Theorem \ref{thm2} which
gives a general bound on its Frobenius error conditionally on bounds
for the stochastic terms.
Theorem~\ref{thm3}, provides bounds on the Frobenius error of our
estimator in closed form. Therefore, we use bounds on the stochastic
terms that we derive in Section \ref{stochastic}. To obtain such
bounds, we use a non-commutative extension of the classical Bernstein
inequality.

In Section \ref{completionunknown}, we consider the case when the
variance of the noise is unknown. Our construction uses the idea of
``square-root'' estimators, first introduced by Belloni \textit{et al.}  \cite
{chernozhukov-square} in the case of the square-root Lasso estimator.
Theorem \ref{thmu3}, shows that our estimator has the same
performances as previously considered estimators which require the
knowledge of the standard deviation of the noise and of the sampling
distribution.

\section{Preliminaries}\label{Preliminaries}
\subsection{Model and sampling scheme}
Let $A_0\in\mathbb{R}^{m_{1}\times m_{2}}$ be an unknown matrix, and
consider the observations $(X_i,Y_i)$ satisfying the trace regression model
%
\begin{equation}
\label{model} Y_{i}=\mathrm{tr}\bigl(X_{i}^{T}A_{0}
\bigr)+\sigma\xi_{i}, \qquad i=1,\ldots,n.
\end{equation}
The noise variables $\xi_{i}$ are independent, with
$\bE(\xi_i)=0$ and $\bE(\xi_i^{2})=1$;
$X_{i}$ are random matrices of dimension $m_{1}\times m_{2}$ and
$\mathrm{tr}(A)$ denotes the trace of the matrix $A$.
Assume that the design matrices $X_i$ are i.i.d. copies of a random
matrix $X$ having distribution $\Pi$ on the set
%
\begin{equation}
\label{basisUSR} \mathcal{X} = \bigl\{e_j(m_1)e_k^{T}(m_2),1
\leq j\leq m_1, 1\leq k\leq m_2 \bigr\},
\end{equation}
where $e_l(m)$ are the canonical basis vectors in $\bR^{m}$. Then, the
problem of estimating $A_0$ coincides with the problem of matrix
completion with random sampling distribution $\Pi$.

One of the particular settings of this problem is the Uniform Sampling
at Random (USR) matrix completion which corresponds to the uniform
distribution $\Pi$. We consider a more general weighted sampling
model. More precisely, let $\pi_{jk}=\mathbb{P}
(X=e_j(m_1)e_k^{T}(m_2) )$ be the probability to observe the
$(j,k)$th entry.
Let us denote by $C_k=\sum_{j=1}^{m_1}\pi_{jk}$
the probability to observe an element from the $k$th column and by
$R_j=\sum_{k=1}^{m_2}\pi_{jk}$ the probability to
observe an element from the $j$th row. Observe that
$\max_{i,j} (C_i,R_j )\geq1/\min(m_1,m_2)$.

As it was shown in \cite{serbro-collaborative}, the trace-norm
penalization fails in the specific situation when the row/column
marginal distribution is such that some columns or rows are sampled
with very high probability (for more details, see \cite
{serbro-collaborative,serbro-learning}). To avoid such a situation, we
need the following assumption on the sampling distribution:
%
\begin{assumption}\label{L}
There exists a positive constant $L\geq1$ such that
\[
\max_{i,j} (C_i,R_j )\leq L/
\min(m_1,m_2).
\]
\end{assumption}
In order to get bounds in the Frobenius norm, we suppose that each
element is sampled with positive probability:
%
\begin{assumption}\label{assPi}
There exists a positive constant $\mu\geq1$ such that
\[
\pi_{jk}\geq(\mu m_1 m_2)^{-1}.
\]
\end{assumption}
In the case of uniform distribution $L=\mu=1$. Let us set $ \Vert
A\Vert_{L_2(\Pi)}^{2}=\bE (\langle A,X\rangle^{2} )$.
Assumption~\ref{assPi} implies that
%
\begin{equation}
\label{ass1} \Vert A\Vert^{2}_{L_2(\Pi)}\geq(m_1m_2
\mu)^{-1}\Vert A\Vert^{2}_{2}.
\end{equation}
%
\subsection{Notation}
We provide a brief summary of the notation used throughout this paper.
Let $A,B$ be matrices in $\mathbb{R}^{m_{1}\times m_{2}}$.
\begin{itemize}
\item We define the \textit{scalar product}
$\langle A,B\rangle=\mathrm{tr}(A^{T}B)$.
\item For $0<q<\infty$ the \textit{Schatten-q} (\textit{quasi-})\textit{norm} of the
matrix $A$ is defined by
\[
\Vert A\Vert_q= \Biggl(\sum_{j=1}^{\min(m_1,m_2)}
\sigma_j(A)^{q} \Biggr)^{1/q} \quad\mbox{and}\quad \Vert A
\Vert=\sigma_1(A),
\]
where $(\sigma_j(A))_j$ are the singular values of $A$ ordered
decreasingly.
\item$\llVert
A\rrVert_{\infty}=\max_{i,j}| a_{ij}|$ where
$A=(a_{ij})$.
\item Let $\pi_{i,j}=\mathbb{P}
(X=e_i(m_1)e_j^{T}(m_2) )$ be the probability to observe the
$(i,j)$th element.
\item For $j=1,\ldots, m_2$,
$C_j=\sum_{i=1}^{m_1}\pi_{ij}$ and for $i=1,\ldots, m_1$,
$R_i=\sum_{j=1}^{m_2}\pi_{ij}$.
\item
$R=\mathrm{diag}(R_1,\ldots,R_{m_1})$ and
$C=\mathrm{diag}(C_1,\ldots,C_{m_2})$.
\item Let $M=\max(m_1,m_2)$,
$m=\min(m_1,m_2)$ and $d=m_1+m_2$.
\item$ \Vert A\Vert_{L_2(\Pi)}^{2}=\bE (\langle A,X\rangle^{2}
)$.\vspace*{1pt}
\item Let
$\{\varepsilon_i\}_{i=1}^{n}$ be an i.i.d. Rademacher sequence and we
define
%
\begin{equation}
\label{stoch1} \Sigma_R=\frac{1}{n} \sum
_{i=1}^{n}\varepsilon_i X_i
\quad\mbox{and}\quad \Sigma=\frac{\sigma}{n}\sum_{i=1}^{n}
\xi_iX_i.
\end{equation}
\item Define the observation operator $\Omega\dvtx \bR^{m_1\times
m_2}\rightarrow\bR^{n} $ as $ (\Omega(A) )_i=
\langle X_i,A \rangle$.
\item$Q(A)=\sqrt{\frac{1}{n}\sum_{i=1}^{n} (Y_i- \langle
X_i,A \rangle )^{2}}$.
\end{itemize}


\section{Matrix completion with known variance of the noise}
\label{completionknown}
In this section, we consider the matrix completion problem when the
variance of the noise is known.
We define the following estimator of $A_0$:
%
\begin{equation}
\label{estimator} \hat{A}=\mathop{\arg\min}_{
\llVert  A\rrVert _{\infty}\leq\mathbf{a}} \Biggl\{
\frac
{1}{n}\sum_{i=1}^{n}
\bigl(Y_i- \langle X_i,A \rangle \bigr)^{2}+
\lambda\Vert A\Vert_1 \Biggr\},
\end{equation}
where $\lambda>0$ is a regularization parameter and $\mathbf{a}$ is
an upper bound on $\llVert  A_0\rrVert_{\infty}$. This is a
restricted version of the matrix LASSO estimator. The matrix LASSO
estimator is based on a trade-off between fitting the target matrix to
the data using least squares and minimizing the nuclear norm and it has
been studied by a number of authors (see, e.g., \cite{candes-plan-noise,rhode-tsybakov-estimation,wainwright-estimation}).

A restricted version of a slightly different estimator, penalised by a
weighted nuclear norm $ \llVert \sqrt{R}A\sqrt{C}\rrVert_1$, was first considered by Negahban and Wainwright in \cite
{wainwright-weighted}. Here $R$ and $C$ are diagonal matrices with
diagonal entries $ \{R_j, j=1,\ldots, m_1 \}$ and $ \{
C_k, k=1,\ldots, m_2 \}$, respectively.
In \cite{wainwright-weighted}, the domain of optimization is the
following one
%
\begin{equation}
\label{domainwainwrite} \biggl\{A\dvtx  \llVert A\rrVert_{\omega(\infty)}
\leq\frac
{\alpha^{*}}{\sqrt{m_1 m_2}} \biggr\},
\end{equation}
where $\alpha^{*}$ is a bound on the ``spikiness ratio'' $\alpha_{sp}=\frac{ \sqrt{m_1m_2}\llVert  A_0\rrVert_{\omega
(\infty)}}{\llVert  A_0\rrVert_{\omega(2)}}$ of the unknown
matrix~$A_0$. Here $\llVert  A\rrVert_{\omega(\infty)}=\llVert \sqrt{R}A\sqrt{C}\rrVert_{\infty}$ and $\llVert
A\rrVert_{\omega(2)}=\llVert \sqrt{R}A\sqrt{C}\rrVert_{2}$. In the particular setting of the uniform sampling (\ref
{domainwainwrite}) gives
\[
\bigl\{A\dvtx  \llVert A\rrVert_{\infty}\leq\alpha \bigr\},
\]
where $\alpha$ is an upper bound on the ``spikiness ratio'' $\frac
{\sqrt{m_1 m_2}\llVert  A_0\rrVert_{\infty}}{\llVert
A_0\rrVert_{2}}$.\vadjust{\goodbreak}

The following theorem gives a general upper bound on the prediction
error of estimator $\hat A$ given by (\ref{estimator}). Its proof is
given in Appendix \ref{proof-thm2}. The stochastic terms $ \llVert
\Sigma\rrVert $ and $ \llVert \Sigma_R\rrVert $ play a
key role in what follows.
%
\begin{theorem}\label{thm2}
Let $X_i$ be i.i.d. with distribution $\Pi$ on $\mathcal{X}$ which
satisfies Assumptions \ref{L} and \ref{assPi} and $\lambda>3 \llVert \Sigma\rrVert $. Assume that $ \llVert  A_0\rrVert_{\infty}\leq\mathbf{a}$ for some constant $\mathbf{a}$. Then,
there exist numerical constants $(c_1,c_2)$ such that
\[
\frac{\Vert\hat A-A_0\Vert_2^{2}}{m_1m_2} \leq\max \biggl\{c_1 \mu^{2}
m_1m_2\rank(A_0) \bigl(
\lambda^{2}+\mathbf{a}^{2} \bigl(\bE \bigl( \llVert
\Sigma_R\rrVert \bigr) \bigr)^{2} \bigr),c_2
\mathbf{a}^{2} \mu\sqrt{\frac{\log(d)}{n}} \biggr\}
\]
with probability at least $1-\frac{2}{d}$, where $d=m_1+m_2$.
\end{theorem}
In order to get a bound in a closed form, we need to obtain suitable
upper bounds on $\bE ( \llVert \Sigma_R\rrVert )$ and, with probability
close to $1$, on $\llVert \Sigma\rrVert $. We will obtain such bounds
in the case of \textit{sub-exponential noise}, that is, under the
following assumption:
%
\begin{assumption}\label{noise}
\[
\max_{i=1,\ldots,n}\bE\exp \bigl(\vert\xi_i\vert/K \bigr)< \infty.
\]
\end{assumption}
Let $K>0$ be a constant
such that $\max_{i=1,\ldots,n}\bE\exp (\vert\xi_i\vert
/K )\leq e$.
The following two lemmas give bounds on $\llVert \Sigma\rrVert $ and $\bE ( \llVert \Sigma_R\rrVert  )$.
We prove them in Section \ref{stochastic} using the non-commutative
Bernstein inequality.
%
\begin{lemma}\label{delta}
Let $X_i$ be i.i.d. with distribution $\Pi$ on $\mathcal{X}$ which
satisfies Assumptions \ref{L} and~\ref{assPi}. Assume that $(\zeta_i)_{i=1}^{n}$ are independent with $\bE(\zeta_i)=0$, $\bE
(\zeta_i^{2} )=1$ and satisfy Assumption~\ref{noise}. Then,
there exists an absolute constant $C^{*}>0$ that depends only on $K$
and such that, for all $t>0$ with probability at least $1-e^{-t}$ we have
%
\begin{equation}
\label{Max} \Biggl\llVert \frac{1}{n}\sum_{i=1}^{n}
\zeta_iX_i\Biggr\rrVert \leq C^{*}\max \biggl
\{\sqrt{\frac{L(t+\log(d))}{mn}},\frac{\log(m)
(t+\log(d) )}{n} \biggr\},
\end{equation}
where $d=m_1+m_2$.
\end{lemma}
%
\begin{lemma}\label{Edelta}
Let $X_i$ be i.i.d. with distribution $\Pi$ on $\mathcal{X}$ which
satisfies Assumptions \ref{L} and~\ref{assPi}. Assume that $(\zeta_i)_{i=1}^{n}$ are independent with $\bE(\zeta_i)=0$, $\bE
(\zeta_i^{2} )=1$ and satisfy Assumption~\ref{noise}. Then, for
$n\geq m\log^{3}(d)/L$, there exists an absolute constant $C^{*}>0$
such that
\[
\bE\Biggl\llVert \frac{1}{n}\sum_{i=1}^{n}
\zeta_iX_i\Biggr\rrVert \leq C^{*}\sqrt{
\frac{2eL\log(d)}{nm}},
\]
where $d=m_1+m_2$.
\end{lemma}
An optimal choice of the parameter $t$ in these lemmas is $t=\log(d)$.
Larger $t$ leads to a slower rate of convergence and a smaller $t$ does
not improve the rate but
makes the concentration probability smaller. With this choice of $t$
the second terms in the maximum in (\ref{Max}) is negligible for
$n>n^{*}$ where $n^{*}=2\log^{2}(d)m/L$. Then, we can choose
%
\begin{equation}
\label{lambda} \lambda=3C^{*}\sigma\sqrt{\frac{2L\log(d)}{mn}},
\end{equation}
where $C^{*}$ is an absolute numerical
constant which depends only on $K$. If $\xi_i$ are $N(0,1)$, then we
can take $C^{*}=6.5$ (see Lemma 4 in \cite{klopp-variance}).
With this choice of $\lambda$, we obtain the following theorem.
%
\begin{theorem}\label{thm3}
Let $X_i$ be i.i.d. with distribution $\Pi$ on $\mathcal{X}$ which
satisfies Assumptions \ref{L} and \ref{assPi}. Assume that $ \llVert  A_0\rrVert_{\infty}\leq\mathbf{a}$ for some constant
$\mathbf{a}$ and that Assumption \ref{noise} holds. Consider the
regularization parameter $\lambda$ satisfying (\ref{lambda}). Then,
there exist a numerical constant $c'$, that depends only on $K$, such that
%
\begin{equation}
\label{revisionthm} \frac{\Vert\hat A-A_0\Vert_2^{2}}{m_1m_2}\leq c' \max
\biggl\{\max \bigl(\sigma^{2},\mathbf{a}^{2}\bigr)
\mu^{2} L\frac{\log(d)\rank
(A_0)M}{n},\mathbf{a}^{2}\mu\sqrt{
\frac{\log(d)}{n}} \biggr\}
\end{equation}
with probability greater than $1-3/d$.
\end{theorem}

\begin{Remarks*}
\textit{Comparison to other works}: An important feature of
our estimator is that its construction requires only an upper bound on
the maximum absolute value of the entries of $A_0$ (and an upper bound
on the variance of the noise). This condition is very mild. Let us
compare this matrix condition and the bound we obtain with some of the
previous works on noisy matrix completion.

We will start with the paper of Keshavan \textit{et al.}  \cite
{keshavan-montanari-matrix}. Their method requires a priori information
on the rank of the unknown matrix as well as a matrix incoherence
assumption (which is stated in terms of the singular vectors of $A_0$).
Under a sampling scheme different from ours (uniform sampling without
replacement) and sub-Gaussian errors, the estimator proposed in \cite
{keshavan-montanari-matrix} satisfies, with high probability, the
following bound
%
\begin{equation}
\label{keshavan} \frac{\Vert\hat A -A_0\Vert_2^{2}}{m_1m_2}\lesssim k^{4}\sqrt {\alpha}
\frac{M}{n}\rank(A_0)\log n.
\end{equation}
The symbol $\lesssim$ means that the inequality holds up to
multiplicative numerical constants, $k=\sigma_{\max}(A_0)/\sigma_{\min}(A_0)$ is the condition number and $\alpha=(m_1\vee
m_2)/(m_1\wedge m_2)$ is the aspect ratio. Comparing (\ref{keshavan})
and (\ref{revisionthm}), we see that our bound is better: it does not
involve the multiplicative coefficient $k^{4}\sqrt{\alpha}$ which can
be big.

Wainwright \textit{et al.}  in \cite{wainwright-weighted} propose an estimator
which uses a priori information on the ``spikiness ratio'' $\alpha_{sp}=\frac{ \sqrt{m_1m_2}\llVert  A_0\rrVert_{\infty
}}{\llVert  A_0\rrVert_2}$ of $A_0$. This method requires
$\alpha_{sp}$ bounded by a constant, say~$\alpha_{*}$, in which case
the estimator proposed in \cite{wainwright-weighted} satisfies the
following bound
%
\begin{equation}
\label{wainwright} \frac{\Vert\hat A -A_0\Vert_{\omega(2)}^{2}}{m_1m_2}\lesssim \alpha_{*}^{2}
\frac{M}{n}\rank(A_0)\log m.
\end{equation}
In the case of uniform sampling and bounded ``spikiness ratio'' this
bound coincides with the bound given by Theorem \ref{thm3}. An
important advantage of our method is that the sampling distribution
does not need to be equal to a product distribution (i.e., $\pi_{ij}$
need not be equal to $R_i C_j$) as is required in \cite{wainwright-weighted}.

The methods proposed in \cite
{Koltchinskii-Tsybakov,klopp-rank,klopp-variance} use the uniform
sampling. Similarly to our construction, an a priori bound on $ \llVert  A_0\rrVert_\infty$ is required. An important difference is
that, in these papers, the bound on $ \llVert  A_0\rrVert_\infty$ is used in the choice of the regularization parameter
$\lambda$. This implies that the convex functional which is minimized
in order to obtain $\hat A$ depends on $\mathbf{a}$. A too large bound
may jeopardize the exactness of the estimation. In our construction,
$\mathbf{a}$ determines the ball over which we are minimizing our
convex functional, which itself is independent of $\mathbf{a}$. Our
estimator achieves the same bound as the estimators proposed in these papers.

\textit{Minimax optimality}: If we consider the matrix
completion setting (i.e., $n\leq m_1 m_2$), then, the maximum in
(\ref{revisionthm}) is given by its first therm. In the case of
Gaussian errors and under the additional assumption that $\pi_{jk}\leq
\frac{\mu_1}{m_1m_2}$ for some\vspace*{1pt} constant $\mu_1\geq1$ this rate of
convergence is minimax optimal (cf. Theorem 5 of \cite
{Koltchinskii-Tsybakov}). This optimality holds for the class of
matrices $\mathcal{A}(r,a)$ defined as follows: for given $r$ and $a$
$A_0\in\mathcal{A}(r,a)$ if and only if the rank of $A_0$ is not
larger than $r$ and all the entries of $A_0$ are bounded in absolute
value by $a$.

\textit{Possible extensions}: The techniques developed in
this paper may also be used to analyse weighted trace norm penalty
similar to one used in \cite{wainwright-weighted,serbro-learning}.
\end{Remarks*}

\section{Matrix completion with unknown variance of the noise}
\label{completionunknown}
In this section, we propose a new estimator for the matrix completion
problem in the case when the variance of the noise $\sigma$ is
unknown. Our construction is inspired
by the square-root Lasso estimator proposed in \cite{chernozhukov-square}.
We define the following estimator of $A_0$:
%
\begin{equation}
\label{estimator-unknown} \hat{A}_{\mathrm{SQ}}=\mathop{\arg\min}_{
\llVert  A\rrVert _{\infty}\leq\mathbf{a}}
\Biggl\{\sqrt {\frac{1}{n}\sum_{i=1}^{n}
\bigl(Y_i- \langle X_i,A \rangle \bigr)^{2}}+
\lambda\Vert A\Vert_1 \Biggr\},
\end{equation}
where $\lambda>0$ is a regularization parameter and $\mathbf{a}$ is
an upper bound on $\llVert  A_0\rrVert_{\infty}$. Note that
the first term of this estimator is the square root of the
data-dependent term of the estimator that we considered in Section \ref
{completionknown}. This is similar to the principle used to
define the square-root Lasso estimator for the usual vector regression model.

Let us set $\rho=\frac{1}{16 \mu m_1m_2\rank(A_0)}$. The
following theorem gives a general upper bound on the prediction error
of the estimator $\hat A_{\mathrm{SQ}}$. Its proof is given in Appendix \ref
{proof-thmu1}.
%
\begin{theorem}\label{thmu1}
Let $X_i$ be i.i.d. with distribution $\Pi$ on $\mathcal{X}$ which
satisfies Assumptions \ref{L} and \ref{assPi}. Assume that $ \llVert  A_0\rrVert_{\infty}\leq\mathbf{a}$ for some constant
$\mathbf{a}$ and $\sqrt{\rho}\geq\lambda\geq3 \llVert \Sigma
\rrVert /Q(A_0)$. Then, there exist numerical constants $c_1'$,
that depends only on $K$, such that with probability at least $1-\frac{2}{d}$
\begin{eqnarray*}
\frac{\llVert \hat A_{\mathrm{SQ}}-A_0\rrVert_2^{2}}{m_1m_2}&\leq& c'_1\max\biggl\{
\mu^{2} m_1m_2\rank(A_0) \bigl(
Q^{2}(A_0)\lambda^{2}+\mathbf{a}^{2}
\bigl(\bE \bigl( \llVert \Sigma_R\rrVert \bigr) \bigr)^{2}
\bigr),
\\
&&\hspace*{35pt} \mathbf{a}^{2}\mu\sqrt{\frac{\log(d)}{n}}\biggr\},
\end{eqnarray*}
where $Q(A_0)=\sigma\sqrt{\frac{1 }{n}\sum_{i=1}^{n}\xi^{2}_i}$.
\end{theorem}
In order to get a bound on the prediction risk in a closed form, we use
the bounds on $\llVert \Sigma\rrVert $ and $\bE ( \llVert
\Sigma_R\rrVert  )$ given by Lemmas \ref{delta} and \ref{Edelta} taking
$t=\log(d)$. It remains to bound $Q(A_0)=\sigma \sqrt{\frac{1
}{n}\sum_{i=1}^{n}\xi^{2}_i}$. We consider the case of sub-Gaussian
noise:
%
\begin{assumption}\label{subG} There exists a constant $K$ such that
\[
\bE \bigl[\exp(t\xi_i) \bigr]\leq\exp \bigl(t^{2}/2K
\bigr)
\]
for all $t>0$.
\end{assumption}
Note that condition $\bE\xi_i^{2}=1$ implies that $K\leq1$.
Under Assumption \ref{subG}, $\xi_{i}^{2}$ are sub-exponential random
variables. Then, the Bernstein inequality for sub-exponential random
variables implies that, there exists a numerical constant $c_3$ such
that, with probability at least $1-2\exp\{-c_3 n\}$, one has
%
\begin{equation}
\label{bQ} 3\sigma/2\geq Q(A_0)\geq\sigma/2.
\end{equation}
Using Lemma \ref{delta} and the right-hand side of (\ref{bQ}), for
$n\geq2\log^{2}(d)m/L$, we can take
%
\begin{equation}
\label{lambdaun} \lambda=6C^{*}\sqrt{\frac{2L\log(d)}{mn}}.
\end{equation}
Note that $\lambda$ \textit{does not depend} on $\sigma$ and
satisfies the two conditions required in Theorem \ref{thmu1}. We have that
%
\begin{equation}
\label{un9} \lambda\geq3 \llVert \Sigma\rrVert /Q(A_0)
\end{equation}
with probability greater then $1-1/d-2\exp\{-c_3 n\}$ and
%
\begin{equation}
\label{un10}
\lambda^{2}\leq\frac{1}{16\mu m_1m_2\rank(A_0)}
\end{equation}
for $n$ large enough, more precisely, for $n$ such that
%
\begin{equation}
\label{nun} n\geq c_4 \mu L M\rank(A_0)\log(d),
\end{equation}
where $c_4=576 (C^{*})^{2}$.
We obtain the following theorem.
%
\begin{theorem}\label{thmu3}
Let $X_i$ be i.i.d. with distribution $\Pi$ on $\mathcal{X}$ which
satisfies Assumptions \ref{L} and \ref{assPi}. Assume that $ \llVert  A_0\rrVert_\infty\leq\mathbf{a}$ for some constant
$\mathbf{a}$ and that Assumption \ref{subG} holds. Consider the
regularization parameter $\lambda$ satisfying (\ref{lambdaun}) and
$n$ satisfying (\ref{nun}). Then, there exist numerical constants
$(c'',c_3)$ such that,
%
\begin{equation}
\frac{\Vert\hat A_{\mathrm{SQ}}-A_0\Vert_2^{2}}{m_1m_2}\leq c''\max \biggl\{ \max\bigl(
\sigma^{2},\mathbf{a}^{2}\bigr)\mu^{2} L
\frac{\log(d)\rank
(A_0)M}{n},\mathbf{a}^{2}\mu\sqrt{\frac{\log(d)}{n}} \biggr\}
\end{equation}
with probability greater than $1-3/d-2\exp\{-c_3 n\}$.
\end{theorem}
Note that condition (\ref{nun}) is not restrictive: indeed the
sampling sizes $n$ satisfying condition (\ref{nun}) are of the same
order of magnitude as those for which the normalized Frobenius error of
our estimator is small. Thus, Theorem \ref{thmu3} shows, that $\hat
A_{\mathrm{SQ}}$ has the same prediction performances as previously proposed
estimators which rely on the knowledge of the standard deviation of the
noise and of the sampling distribution.

\section{Bounds on the stochastic errors}\label{stochastic}
In this section, we will obtain the upper bounds for the stochastic
errors $\llVert \Sigma_R\rrVert $ and $\bE ( \llVert \Sigma_R\rrVert  )$ defined in (\ref{stoch1}). In
order to obtain such bounds, we use the matrix version of Bernstein's
inequality.
The following proposition is obtained by an extension of Theorem 4 in
\cite{koltchinskii-remark} to rectangular matrices via self-adjoint
dilation (cf., for example, 2.6 in \cite{tropp-user}).
Let $Z_1,\ldots,Z_n$ be independent random matrices with dimensions
$m_1\times m_2$. Define
\[
\sigma_Z=\max \Biggl\{\Biggl\llVert \frac{1}{n}\sum
_{i=1}^{n}\bE \bigl(Z_iZ^{T}_i
\bigr)\Biggr\rrVert^{1/2}, \Biggl\llVert \frac{1}{n}\sum
_{i=1}^{n} \bE \bigl(Z_i^{^{T}}Z_i
\bigr)\Biggr\rrVert^{1/2} \Biggr\}
\]
and
\[
U_i=\inf \bigl\{K>0\dvtx  \bE\exp \bigl(\Vert Z_i\Vert/K \bigr)\leq
e \bigr\}.
\]

\begin{proposition}\label{pr1}
Let $Z_1,\ldots,Z_n$ be independent random matrices with dimensions
$m_1\times m_2$ that satisfy $\bE(Z_i)=0$. Suppose that $U_i<U$ for
some constant $U$ and all $i=1,\ldots,n$. Then, there exists an
absolute constant $c^{*}$, such that, for all $t>0$, with probability
at least $1-e^{-t}$ we have
\[
\Biggl\llVert \frac{1}{n}\sum_{i=1}^{n}Z_i
\Biggr\rrVert \leq c^{*}\max \biggl\{ \sigma_Z\sqrt{
\frac{t+\log(d)}{n}},U \biggl(\log\frac{U}{\sigma_Z} \biggr)\frac{t+\log(d)}{n}
\biggr\},
\]
where $d=m_1+m_2$.
\end{proposition}
%
\subsection{\texorpdfstring{Proof of Lemma \protect\ref{delta}}{Proof of Lemma 5}}
We apply Proposition \ref{pr1} to $Z_i=\zeta_iX_i$.
We first estimate $\sigma_Z$ and $U$. Note that $Z_i$ is a zero-mean
random matrix which satisfies
\[
\llVert Z_i\rrVert \leq\vert\zeta_i\vert.
\]
Then, Assumption \ref{noise} implies that there exists a constant $K$
such that $U_i\leq K$ for all $i=1,\ldots,n$.
We compute
\[
\bE \bigl(Z_iZ^{T}_i \bigr)= R \quad\mbox{and}\quad
\bE \bigl(Z^{T}_iZ_i \bigr)= C,
\]
where $C$ (resp., $R$) is the diagonal matrix with $C_k$ (resp., $R_j$)
on the diagonal. 
This and the fact that the $X_i$ are i.i.d. imply that
\[
\sigma_Z^{2}=\max_{i,j} (C_i,R_j
)\leq L/m.
\]
Note that $\max_{i,j} (C_i,R_j )\geq1/m$ which implies
that $\log ( K/\sigma_Z )\leq\log (Km ) $ and
the statement of Lemma \ref{delta} follows.

\subsection{\texorpdfstring{Proof of Lemma \protect\ref{Edelta}}{Proof of Lemma 6}}
The proof follows the lines of the proof of Lemma 7 in \cite
{klopp-rank}. For sake of completeness, we give it here.
Set $t^{*}=\frac{Ln}{m\log^{2}(m)}-\log(d)$. $t^{*}$ is the\vspace*{1pt} value
of $t$ such that the two terms in (\ref{Max}) are equal.
Note that Lemma \ref{delta} implies that
%
\begin{equation}
\label{proba1} \mathbb P \Biggl(\Biggl\llVert \frac{1}{n}\sum
_{i=1}^{n}\zeta_iX_i
\Biggr\rrVert > t \Biggr)\leq d\exp\bigl\{-t^{2} nm/ \bigl(
\bigl(C^{*}\bigr)^{2}L \bigr)\bigr\} \qquad\mbox{for } t\leq
t^{*}
\end{equation}
and
%
\begin{equation}
\label{proba2} \mathbb P \Biggl(\Biggl\llVert \frac{1}{n}\sum
_{i=1}^{n}\zeta_iX_i
\Biggr\rrVert > t \Biggr)\leq d\exp\bigl\{-t n/\bigl(C^{*}\log( m)
\bigr)\bigr\} \qquad\mbox{for } t\geq t^{*}.
\end{equation}
We set $\nu_1=nm/ ((C^{*})^{2}L )$, $\nu_2=n/(C^{*}\log
(m))$. By H{\"o}lder's inequality, we get
\[
\bE\Biggl\llVert \frac{1}{n}\sum_{i=1}^{n}
\zeta_iX_i\Biggr\rrVert \leq \Biggl(\bE\Biggl\llVert
\frac{1}{n}\sum_{i=1}^{n}
\zeta_iX_i\Biggr\rrVert^{2\log(
d)}
\Biggr)^{1/(2\log(d))}.
\]
The inequalities (\ref{proba1}) and (\ref{proba2}) imply that
%
\begin{eqnarray}
\label{estEM}
&&
\Biggl(\bE\Biggl\llVert \frac{1}{n}\sum
_{i=1}^{n}\zeta_iX_i\Biggr
\rrVert^{2\log
(d)} \Biggr)^{1/2\log(d)} \nonumber\\
&&\quad= \Biggl( \int^{+\infty}_{0}
\mathbb P \Biggl(\Biggl\llVert \frac{1}{n}\sum_{i=1}^{n}
\zeta_iX_i\Biggr\rrVert > t^{1/(2\log
(d))} \Biggr)
\,\mathrm{d}t \Biggr)^{1/2\log(d)}
\nonumber\\[-8pt]\\[-8pt]
&&\quad\leq \biggl(d \int^{+\infty}_{0}\exp\bigl
\{-t^{1/\log
(d)}\nu_1\bigr\}\,\mathrm{d}t+d \int
^{+\infty}_{0}\exp\bigl\{-t^{1/(2\log
(d)}
\nu_2\bigr\}\,\mathrm{d}t \biggr)^{1/2\log(d)}
\nonumber\\
&&\quad\leq \sqrt{e} \bigl(\log(d)\nu_1^{-\log(d)}\Gamma
\bigl(\log(d)\bigr)+2\log(d) \nu_2^{-2\log(d)}\Gamma\bigl(2\log(d)
\bigr) \bigr)^{1/(2\log(d))}.
\nonumber
\end{eqnarray}
The Gamma-function satisfies the following bound:
%
\begin{equation}
\label{Gamma} \mbox{for } x\geq2\qquad \Gamma(x)\leq \biggl(\frac
{x}{2}
\biggr)^{x-1}
\end{equation}
(see, e.g., \cite{klopp-rank}). Plugging this into (\ref{estEM}), we compute
\begin{eqnarray*}
&&
\bE\Biggl\llVert \frac{1}{n}\sum_{i=1}^{n}
\zeta_iX_i\Biggr\rrVert \\
&&\quad\leq \sqrt {e}\bigl(\bigl(
\log(d)\bigr)^{\log(d)}\nu_1^{-\log(d)}2^{1-\log(d)}
+2\bigl(\log(d)\bigr)^{2\log(d)}\nu_2^{-2\log(d)}
\bigr)^{1/(2\log(d))}.
\end{eqnarray*}
Observe that $n>n^{*}$ implies
$\nu_1\log(d)\leq\nu_2^{2}$ and we obtain
%
\begin{equation}
\label{estEM-1} \bE\Biggl\llVert \frac{1}{n}\sum
_{i=1}^{n}\zeta_iX_i\Biggr
\rrVert \leq\sqrt {\frac{2e\log(d)}{\nu_1}}.
\end{equation}
We conclude the proof by plugging $\nu_1=nm/ ((C^{*})^{2}L
)$ into (\ref{estEM-1}).

\begin{appendix}
\section{\texorpdfstring{Proof of Theorem \protect\ref{thm2}}{Proof of Theorem 3}}\label{proof-thm2}
It follows from the definition of the estimator $\hat A$ that
\[
\frac{1}{n}\sum_{i=1}^{n}
\bigl(Y_i- \langle X_i,\hat A \rangle
\bigr)^{2}+\lambda\Vert\hat A\Vert_1\leq\frac{1}{n}
\sum_{i=1}^{n} \bigl(Y_i-
\langle X_i,A_0 \rangle \bigr)^{2}+\lambda
\Vert A_0\Vert_1,
\]
which, using (\ref{model}), implies
\[
\frac{1}{n}\sum_{i=1}^{n} \bigl(
\langle X_i,A_0 \rangle +\sigma \xi_i-
\langle X_i,\hat A \rangle \bigr)^{2}+\lambda \Vert\hat A
\Vert_1\leq\frac{\sigma^{2}}{n}\sum_{i=1}^{n}
\xi_i^{2}+\lambda \Vert A_0
\Vert_1.
\]
Hence,
\[
\frac{1}{n}\sum_{i=1}^{n} \langle
X_i,A_0-\hat A \rangle^{2}+2 \langle
\Sigma,A_0-\hat A \rangle
+\lambda\Vert \hat A\Vert_1\leq\lambda\Vert A_0
\Vert_1,
\]
where $\Sigma=\frac{\sigma}{n}\sum_{i=1}^{n}\xi_iX_i$. Then, by
the duality
between the nuclear and the operator norms, we obtain
%
\setcounter{equation}{23}
\renewcommand{\theequation}{\arabic{equation}}
\begin{equation}
\label{1} \frac{1}{n} \bigl\llVert \Omega (A_0-\hat A )\bigr
\rrVert_2^{2}+\lambda\Vert\hat A\Vert_1\leq2
\llVert \Sigma\rrVert \Vert A_0-\hat A\Vert_1+\lambda
\Vert A_0\Vert_1.
\end{equation}

Let $P_S$ be the projector on the linear vector subspace $S$ and
let $S^\bot$ be the orthogonal complement of $S$.
Let $u_j(A)$ and $v_j(A)$ denote, respectively, the \textit{left} and
\textit{right} orthonormal \textit{singular vectors} of $A$. $S_1(A)$
is the linear span of $\{u_j(A)\}$, $S_2(A)$ is the linear span of $\{
v_j(A)\}$.
We set
%
\begin{equation}
\label{projector} \mathbf P_A^{\bot}(B)=P_{S_1^{\bot}(A)}BP_{S_2^{\bot}(A)}
\quad\mbox{and}\quad \mathbf P_A(B)=B- \mathbf P_A^{\bot}(B).
\end{equation}

By definition of $\mathbf P_{A_0}^{\bot}$, for any matrix $B$, the
singular vectors of $\mathbf P_{A_0}^{\bot}(B)$ are orthogonal to the
space spanned by the singular vectors of $A_0$. This implies that
$\llVert  A_0+\mathbf P_{A_0}^{\bot}(\hat A-A_0) \rrVert_1=\llVert  A_0 \rrVert_1+\llVert \mathbf P_{A_0}^{\bot
}(\hat A-A_0) \rrVert_1$. Then
%
\begin{eqnarray}
\label{ineq} \llVert \hat A\rrVert_1&=&\llVert A_0 +\hat
A-A_0\rrVert_1
=\bigl\llVert A_0 +\mathbf P_{A_0}^{\bot}(\hat
A-A_0)+\mathbf P_{A_0}(\hat A-A_0)\bigr
\rrVert_1
\nonumber\\
&\geq&\bigl\llVert A_0 +\mathbf P_{A_0}^{\bot}(\hat
A-A_0)\bigr\rrVert_1-\bigl\llVert \mathbf
P_{A_0}(\hat A-A_0)\bigr\rrVert_1
\\
&=&\llVert A_0 \rrVert_1+\bigl\llVert \mathbf
P_{A_0}^{\bot
}(\hat A-A_0)\bigr
\rrVert_1-\bigl\llVert \mathbf P_{A_0}(\hat
A-A_0)\bigr\rrVert_1.
\nonumber
\end{eqnarray}
Note that from (\ref{ineq}), we get
%
\begin{equation}
\label{un2correction} \llVert A_0
\rrVert_1-\llVert \hat A\rrVert_1\leq \bigl\llVert
\mathbf P_{A_0}(A_0-\hat A)\bigr\rrVert_1-\bigl
\llVert \mathbf P_{A_0}^{\bot}(A_0-\hat A)\bigr
\rrVert_1.
\end{equation}
This, the triangle inequality and $\lambda\geq3\llVert \Sigma
\rrVert $ lead to
%
\begin{eqnarray}
\label{2} \frac{1}{n} \bigl\llVert \Omega (A_0-\hat A )\bigr
\rrVert_2^{2}&\leq&2\llVert \Sigma\rrVert \bigl\llVert
\mathbf P_{A_0} (A_0-\hat A )\bigr\rrVert_1+
\lambda\bigl\llVert \mathbf P_{A_0} ( A_0-\hat A )\bigr
\rrVert_1
\nonumber\\[-8pt]\\[-8pt]
&\leq& \frac{5}{3}\lambda\bigl\llVert \mathbf P_{A_0}
(A_0-\hat A )\bigr\rrVert_1.
\nonumber
\end{eqnarray}
Since $\mathbf P_A(B)=P_{S_1^{\bot}(A)}BP_{S_2(A)}+ P_{S_1(A)}B$
and $\rank(P_{S_i(A)}B)\leq\rank(A)$ we have that $\rank(\mathbf
P_A(B))\leq2\rank(A)$. 
From (\ref{2}), we compute
%
\begin{equation}
\label{3} \frac{1}{n} \bigl\llVert \Omega (A_0-\hat A )\bigr
\rrVert_2^{2}\leq\frac{5}{3} \lambda\sqrt{2
\rank(A_0)}\llVert \hat A-A_0 \rrVert_2.
\end{equation}

For a $0<r\leq m$, we consider the following constrain set
%
\begin{equation}
\label{constrain} \mathcal{C}(r)= \biggl\{A\in\mathbb{R}^{m_1\times m_2}\dvtx  \llVert
A\rrVert_{\infty}=1, \llVert A\rrVert_{L_2(\Pi
)}^{2}\geq
\sqrt{\frac{64 \log(d)}{\log (6/5 ) n}}, \llVert A\rrVert_{1}\leq\sqrt{r} \llVert A
\rrVert_{2} \biggr\}.
\end{equation}
Note that the condition $\llVert  A\rrVert_{1}\leq\sqrt{r}
\llVert  A\rrVert_{2}$ is satisfied if $\rank(A)\leq r$.

The following lemma shows that for matrices $A\in\mathcal{C}(r)$ the
observation operator $\Omega$ satisfies some approximative restricted
isometry. Its proof is given in Appendix \ref{proof-thm1}.
%
\begin{lemma}\label{thm1}
Let $X_i$ be i.i.d. with distribution $\Pi$ on $\mathcal{X}$ which
satisfies Assumptions \ref{L} and \ref{assPi}. Then, for all $A\in
\mathcal{C}(r)$
\[
\frac{1}{n}\bigl\llVert \Omega(A)\bigr\rrVert^{2}_{2}
\geq\frac
{1}{2}\Vert A\Vert_{L_2(\Pi)}^{2}-44 \mu
rm_1m_2 \bigl(\bE \bigl( \llVert \Sigma_R
\rrVert \bigr) \bigr)^{2}
\]
with probability at least $1-\frac{2}{d}$.
\end{lemma}

We need the following auxiliary lemma which is proven in Appendix \ref{pl2}.
%
\begin{lemma}\label{l2} If $\lambda>3\llVert \Sigma\rrVert $
\[
\bigl\llVert \mathbf P_{A_0}^{\bot}(\hat A-A_0)
\bigr\rrVert_1\leq 5\bigl\llVert \mathbf P_{A_0}(\hat
A-A_0) \bigr\rrVert_1.
\]
\end{lemma}

Lemma \ref{l2} implies that
%
\begin{equation}
\label{revisioncondition} \llVert \hat A-A_0
\rrVert_{1}\leq6\bigl\llVert \mathbf P_{A_0}(\hat
A-A_0) \bigr\rrVert_1\leq\sqrt{72 \rank(A_0)}
\llVert \hat A-A_0\rrVert_{2}.
\end{equation}
Set $a=\llVert \hat A-A_0\rrVert_{\infty}$. By definition of
$\hat A$, we have that $ a\leq2\mathbf{a}$.
We now consider two cases, depending on whether the matrix $\frac
{1}{a} (\hat A-A_0 )$ belongs to the set $ \mathcal{C}
(72 \rank(A_0) )$ or not.

\textit{Case} 1: Suppose first that $\llVert \hat A-A_0\rrVert_{L_2(\Pi)}^{2} < a^{2}\sqrt{\frac{64 \log(d)}{\log
(6/5 ) n}}$, then (\ref{ass1}) implies that
%
\begin{equation}
\label{o1} \frac{\llVert \hat A-A_0\rrVert^{2}_2}{m_1m_2}\leq4\mathbf {a}^{2} \mu\sqrt{
\frac{64 \log(d)}{\log (6/5 ) n}}
\end{equation}
and we get the statement of Theorem \ref{thm2} in this case.

\textit{Case} 2: It remains to consider the case $\llVert \hat
A-A_0\rrVert_{L_2(\Pi)}^{2} \geq a^{2}\sqrt{\frac{64 \log
(d)}{\log (6/5 ) n}}$. Then (\ref{revisioncondition})
implies that $\frac{1}{a} (\hat A-A_0 )\in\mathcal
{C} (72 \rank(A_0) )$ and we can apply Lemma \ref{thm1}.
From Lemma \ref{thm1} and (\ref{3}), we obtain that
with probability at least $1-\frac{2}{d}$ one has
\begin{eqnarray*}
\tfrac{1}{2}\Vert\hat A-A_0\Vert_{L_2(\Pi)}^{2}
&\leq& \tfrac{5}{3}\lambda\sqrt{2\rank(A_0)}\llVert \hat
A-A_0 \rrVert_2+3168 \mu a^{2}
\rank(A_0)m_1m_2 \bigl(\bE \bigl( \llVert
\Sigma_R\rrVert \bigr) \bigr)^{2}
\\
&\leq& 6\lambda^{2}\mu m_1m_2
\rank(A_0)+\tfrac{1}{4} (m_1m_2\mu
)^{-1}\llVert \hat A-A_0 \rrVert^{2}_2
\\
&&{}+3168 \mu a^{2}\rank(A_0)m_1m_2
\bigl(\bE \bigl( \llVert \Sigma_R\rrVert \bigr)
\bigr)^{2}.
\end{eqnarray*}
Now (\ref{ass1}) and $ a\leq2\mathbf{a}$ imply that, there exist
numerical constants $c_1$ such that
\[
\Vert\hat A-A_0\Vert_2^{2}\leq
c_1 (\mu m_1m_2 )^{2}
\rank(A_0) \bigl( \lambda^{2} +\mathbf{a}^{2}
\bigl(\bE \bigl( \llVert \Sigma_R\rrVert \bigr) \bigr)^{2}
\bigr),
\]
which, together with (\ref{o1}), leads to the statement of the Theorem
\ref{thm2}.
\section{\texorpdfstring{Proof of Lemma \protect\ref{thm1}}{Proof of Lemma 12}}\label{proof-thm1}
The main lines of this proof are close to those of the proof of Theorem
1 in \cite{wainwright-weighted}. Set $\mathcal{E}=44 \mu
rm_1m_2 (\bE ( \llVert \Sigma_R\rrVert
) )^{2}$.
We will show that the probability of the following ``bad'' event is small
\[
\mathcal{B}= \biggl\{\exists A\in\mathcal{C}(r) \mbox{ such that } \biggl\llvert
\frac{1}{n} \bigl\llVert \Omega(A) \bigr\rrVert_{2}^{2}-
\Vert A\Vert_{L_2(\Pi)}^{2}\biggr\rrvert > \frac{1}{2}\Vert A
\Vert_{L_2(\Pi)}^{2}+ \mathcal{E} \biggr\}.
\]
Note that $\mathcal{B}$ contains the complement of the event that we
are interested in.

In order to estimate the probability of $\mathcal{B,}$ we use a
standard peeling argument. Let $\nu=\sqrt{\frac{64 \log(d)}{\log
(6/5 ) n}}$ and $\alpha=\frac{6}{5}$. For $l\in\mathbb
N$ set
\[
S_l= \bigl\{A\in\mathcal{C}(r)\dvtx  \alpha^{l-1}\nu\leq\Vert
A\Vert_{L_2(\Pi)}^{2}\leq\alpha^{l}\nu \bigr\}.
\]
If the event $\mathcal{B}$ holds for some matrix $A\in\mathcal
{C}(r)$, then $A$ belongs to some $S_l$ and
%
\setcounter{equation}{32}
\renewcommand{\theequation}{\arabic{equation}}
\begin{eqnarray}
\label{Bl} \biggl\llvert \frac{1}{n} \bigl\llVert \Omega(A) \bigr
\rrVert_{2}^{2}-\Vert A\Vert_{L_2(\Pi)}^{2}
\biggr\rrvert &>& \frac{1}{2}\Vert A\Vert_{L_2(\Pi)}^{2}+
\mathcal{E}
\nonumber
\\
&>& \frac{1}{2}\alpha^{l-1}\nu+ \mathcal{E}
\\
&=& \frac{5}{12}\alpha^{l}\nu+ \mathcal{E}.
\nonumber
\end{eqnarray}
For each $T>\nu$ consider the following set of matrices
\[
\mathcal{C}(r,T)= \bigl\{A\in\mathcal{C}(r)\dvtx  \llVert A\rrVert_{L_2(\Pi)}^{2}
\leq T \bigr\}
\]
and the following event
\[
\mathcal{B}_l= \biggl\{\exists A\in\mathcal{C}\bigl(r,
\alpha^{l}\nu\bigr)\dvtx  \biggl\llvert \frac{1}{n} \bigl\llVert
\Omega(A) \bigr\rrVert_{2}^{2}-\Vert A\Vert_{L_2(\Pi)}^{2}
\biggr\rrvert > \frac
{5}{12}\alpha^{l}\nu+ \mathcal{E} \biggr\}.
\]
Note that $A\in S_l$ implies that $A\in\mathcal{C}(r,\alpha^{l}\nu
)$. Then (\ref{Bl}) implies that $\mathcal{B}_l$ holds and we get
$\mathcal{B}\subset\bigcup\mathcal{B}_l$. Thus, it is enough to
estimate the probability of the simpler event $\mathcal{B}_l$ and then
apply the union bound. Such an estimation is given by the following
lemma. Its proof is given in Appendix \ref{pl1}. Let
\[
Z_T=\sup_{A\in\mathcal{C}(r,T)}\biggl\llvert \frac{1}{n} \bigl\llVert
\Omega(A) \bigr\rrVert_{2}^{2}-\Vert A\Vert_{L_2(\Pi)}^{2}
\biggr\rrvert.
\]

\begin{lemma}\label{l1}
Let $X_i$ be i.i.d. with distribution $\Pi$ on $\mathcal{X}$ which
satisfies Assumptions \ref{L} and \ref{assPi}. Then,
\[
\mathbb P \bigl(Z_T>\tfrac{5}{12}T+ 44 \mu
rm_1m_2 \bigl(\bE \bigl( \llVert \Sigma_R
\rrVert \bigr) \bigr)^{2} \bigr)\leq \exp\bigl(-c_5 n
T^{2}\bigr),
\]
where $c_5=\frac{1}{128}$.
\end{lemma}
Lemma \ref{l1} implies that $\mathbb P (\mathcal{B}_l
)\leq\exp(-c_5 n \alpha^{2l}\nu^{2})$. Using the union bound, we obtain
\begin{eqnarray*}
\mathbb P (\mathcal{B} )&\leq& \sum_{l=1}^{\infty}
\mathbb P (\mathcal{B}_l )
\\
&\leq& \sum_{l=1}^{\infty}\exp
\bigl(-c_5 n \alpha^{2l} \nu^{2}\bigr)
\\
&\leq& \sum_{l=1}^{\infty}\exp \bigl(- \bigl(2
c_5 n \log(\alpha) \nu^{2} \bigr)l \bigr),
\end{eqnarray*}
where we used $e^{x}\geq x$. We finally compute for $\nu=\sqrt{\frac
{64 \log(d)}{\log (6/5 ) n}}$
\[
\mathbb P (\mathcal{B} )\leq\frac{\exp (-2 c_5
n \log(\alpha) \nu^{2} )}{1-\exp (-2 c_5 n \log
(\alpha) \nu^{2} )}=\frac{\exp (-\log(d)
)}{1-\exp (-\log(d) )}.
\]
This completes the proof of Lemma \ref{thm1}.
\begin{Remark*}
As we mentioned in the beginning, the main lines of
this proof are close to those of the proof of Theorem 1 in \cite
{wainwright-weighted}. Let us briefly discuss the main differences
between these two proofs. Similarly to Theorem 1 in \cite
{wainwright-weighted} we prove a kind of ``restricted strong
convexity'' on a constrain set. However, our constrain set defined by
(\ref{constrain}) is quite different from the one introduced in \cite
{wainwright-weighted}: 
\[
\mathcal{C}(n;c_0)= \biggl\{A\in\mathbb{R}^{m_1\times m_2}\dvtx
\frac{\sqrt{m_1m_2} \llVert  A \rrVert_1 \llVert
A\rrVert_\infty}{ \llVert  A\rrVert^{2}_2}\leq\frac
{1}{c_0}\sqrt{\frac{n}{d \log(d)}} \biggr\}.
\]
The present proof is also less involved (e.g., we do not need
use the covering argument used in \cite{wainwright-weighted}). One
important ingredient of our proof is a more efficient control of $\bE
\llVert \Sigma_R\rrVert $ given by Lemma \ref{Edelta}
(compare with Lemma 6 in \cite{wainwright-weighted}).
\end{Remark*}

\section{\texorpdfstring{Proof of Lemma \protect\ref{l1}}{Proof of Lemma 14}}\label{pl1}

Our approach is standard: first we show that $Z_T$ concentrates around
its expectation and then we upper bound the expectation.

By definition, $Z_T=\sup_{A\in\mathcal{C}(r,T)}\llvert \frac
{1}{n} \sum_{i=1}^{n} \langle X_i,A \rangle^{2}-\bE
(
\langle X,A \rangle^{2} )\rrvert $. Massart's
concentration inequality (see, e.g., \cite{sara}, Theorem 14.2)
implies that
%
\setcounter{equation}{33}
\renewcommand{\theequation}{\arabic{equation}}
\begin{equation}
\label{concentration} \mathbb P \bigl(Z_T\geq\bE ( Z_T
)+\tfrac{1}{9} \bigl(\tfrac{5}{12}T \bigr) \bigr)\leq\exp
\bigl(-c_5 n T^{2} \bigr),
\end{equation}
where $c_5=\frac{1}{128}$.
Next, we bound the expectation $\bE ( Z_T )$. Using a
standard symmetrization argument (see, e.g., \cite{koltchiskii-stflour}, Theorem
2.1), we obtain
\begin{eqnarray*}
\bE ( Z_T )&=& \bE \Biggl(\sup_{A\in
\mathcal{C}(r,T)}\Biggl\llvert
\frac{1}{n} \sum_{i=1}^{n} \langle
X_i,A \rangle^{2}-\bE \bigl( \langle X,A
\rangle^{2} \bigr)\Biggr\rrvert \Biggr)
\\
&\leq& 2\bE \Biggl(\sup_{A\in
\mathcal{C}(r,T)}\Biggl\llvert \frac{1}{n} \sum
_{i=1}^{n} \varepsilon_i \langle
X_i,A \rangle^{2} \Biggr\rrvert \Biggr),
\end{eqnarray*}
where $\{\varepsilon_i\}_{i=1}^{n}$ is an i.i.d. Rademacher sequence. The
assumption $\llVert  A\rrVert_{\infty}=1$ implies $\llvert  \langle X_i,A \rangle\rrvert \leq1$. Then, the
contraction inequality (see, e.g., \cite{koltchiskii-stflour}) yields
\[
\bE ( Z_T )\leq8\bE \Biggl(\sup_{A\in
\mathcal{C}(r,T)}\Biggl\llvert
\frac{1}{n} \sum_{i=1}^{n}
\varepsilon_i \langle X_i,A \rangle\Biggr\rrvert \Biggr)=8
\bE \Bigl(\sup_{A\in\mathcal{C}(r,T)}\bigl\llvert \langle \Sigma_R,A \rangle
\bigr\rrvert \Bigr),
\]
where $\Sigma_R=\frac{1}{n} \sum_{i=1}^{n}\varepsilon_i X_i$. For
$A\in
\mathcal{C}(r,T)$, we have that
\begin{eqnarray*}
\llVert A \rrVert_1 &\leq& \sqrt{r}\llVert A\rrVert_2
\\
&\leq& \sqrt{\mu rm_1m_2}\llVert A\rrVert_{L_2(\Pi)}
\\
&\leq& \sqrt{\mu m_1m_2 r T},
\end{eqnarray*}
where we have used (\ref{ass1}). Then, by the duality between nuclear
and operator norms, we compute
\[
\bE ( Z_T )\leq8\bE \Bigl(\sup_{\llVert  A
\rrVert _1\leq\sqrt{\mu m_1m_2rT}}\bigl\llvert \langle
\Sigma_R,A \rangle\bigr\rrvert \Bigr)\leq8\sqrt{\mu
m_1m_2rT} \bE \bigl( \llVert \Sigma_R
\rrVert \bigr).
\]
Finally, using
\[
\tfrac{1}{9} \bigl(\tfrac{5}{12}T \bigr)+8\sqrt{\mu
m_1m_2rT} \bE \bigl( \llVert \Sigma_R
\rrVert \bigr)\leq \bigl(\tfrac{1}{9}+\tfrac{8}{9} \bigr)
\tfrac{5}{12}T+44\mu rm_1m_2 \bigl(\bE \bigl(
\llVert \Sigma_R\rrVert \bigr) \bigr)^{2}
\]
and the concentration bound (\ref{concentration}), we obtain that
\[
\mathbb P \bigl(Z_T>\tfrac{5}{12}T+ 44 \mu
rm_1m_2 \bigl(\bE \bigl( \llVert \Sigma_R
\rrVert \bigr) \bigr)^{2} \bigr)\leq \exp\bigl(-c_5 n
T^{2}\bigr)
\]
with $c_5=\frac{1}{128}$ as stated.

\section{\texorpdfstring{Proof of Theorem \protect\ref{thmu1}}{Proof of Theorem 8}}\label{proof-thmu1}
Let us set $\Delta=A_0-\hat A_{\mathrm{SQ}}$ and $Q(A)=\sqrt{\frac{1}{n}\sum_{i=1}^{n}
(Y_i- \langle X_i,A \rangle )^{2}}$. We have that
\begin{eqnarray*}
Q^{2}(\hat A_{\mathrm{SQ}})-Q^{2}(A_0)&=&
\frac{1}{n} \bigl\llVert \Omega(\Delta )\bigr\rrVert_2^{2}+2
\Biggl\langle\frac{\sigma}{n}\sum_{i=1}^{n}
\xi_iX_i,\Delta \Biggr\rangle
\\
&=& \frac{1}{n} \bigl\llVert \Omega (\Delta )\bigr\rrVert_2^{2}+2
\langle\Sigma,\Delta \rangle,
\end{eqnarray*}
where $\Sigma=\frac{\sigma}{n}\sum_{i=1}^{n}\xi_iX_i$. This implies
%
\setcounter{equation}{34}
\renewcommand{\theequation}{\arabic{equation}}
\begin{equation}
\label{un1} \frac{1}{n} \bigl\llVert \Omega (\Delta )\bigr
\rrVert_2^{2}=-2 \langle\Sigma,\Delta \rangle+ \bigl(Q(\hat
A_{\mathrm{SQ}})-Q(A_0) \bigr) \bigl(Q(\hat A_{\mathrm{SQ}})+Q(A_0)
\bigr).
\end{equation}
We need the following auxiliary lemma which is proven in Appendix
\ref{plu1} ($\mathbf P_{A_0}^{\bot}$ and $\mathbf P_{A_0}$ are
defined in (\ref{projector})).
%
\begin{lemma}\label{lu1} If $\lambda>3 \llVert \Sigma\rrVert /Q(A_0)$, then
\[
\bigl\llVert \mathbf P_{A_0}^{\bot}(\Delta) \bigr
\rrVert_1\leq2\bigl\llVert \mathbf P_{A_0}(\Delta) \bigr
\rrVert_1,
\]
where $\Delta=\hat A_{\mathrm{SQ}}-A_0$.
\end{lemma}
Note that from (\ref{ineq}) we get
%
\begin{equation}
\label{un2} \llVert A_0 \rrVert_1-\llVert \hat
A_{\mathrm{SQ}}\rrVert_1\leq\bigl\llVert \mathbf P_{A_0}(
\Delta)\bigr\rrVert_1-\bigl\llVert \mathbf P_{A_0}^{\bot}(
\Delta)\bigr\rrVert_1.
\end{equation}
The definition of $\hat A_{\mathrm{SQ}}$ and (\ref{un2}) imply that
%
\begin{eqnarray}
\label{un3} Q(A_0)+Q(\hat A_{\mathrm{SQ}})&\leq&2Q(A_0)+
\lambda \bigl(\llVert A_0 \rrVert_1-\llVert \hat
A_{\mathrm{SQ}}\rrVert_1 \bigr)
\nonumber\\[-8pt]\\[-8pt]
&\leq& 2Q(A_0)+\lambda \bigl(\bigl\llVert \mathbf P_{A_0}(
\Delta)\bigr\rrVert_1-\bigl\llVert \mathbf P_{A_0}^{\bot}(
\Delta)\bigr\rrVert_1 \bigr)
\nonumber
\end{eqnarray}
and
%
\begin{eqnarray}
\label{un4} Q(\hat A_{\mathrm{SQ}})-Q(A_0)&\leq&\lambda \bigl(
\llVert A_0 \rrVert_1-\llVert \hat A_{\mathrm{SQ}}
\rrVert_1 \bigr)
\nonumber
\\
&\leq& \lambda \bigl(\bigl\llVert \mathbf P_{A_0}(\Delta)\bigr
\rrVert_1-\bigl\llVert \mathbf P_{A_0}^{\bot}(
\Delta)\bigr\rrVert_1 \bigr)
\\
&\leq&\lambda \bigl(2\bigl\llVert \mathbf P_{A_0}(\Delta)\bigr
\rrVert_1-\bigl\llVert \mathbf P_{A_0}^{\bot}(
\Delta)\bigr\rrVert_1 \bigr).
\nonumber
\end{eqnarray}
Lemma \ref{lu1} implies that $2\llVert \mathbf P_{A_0}(\Delta
)\rrVert_1-\llVert \mathbf P_{A_0}^{\bot}(\Delta)\rrVert_1\geq0$.
From (\ref{un3}) and (\ref{un4}), we compute
%
\begin{eqnarray}
\label{un5} && \bigl(Q(\hat A_{\mathrm{SQ}})-Q(A_0) \bigr) \bigl(Q(
\hat A_{\mathrm{SQ}})+Q(A_0) \bigr)
\nonumber\\
&&\quad\leq\lambda \bigl(2\bigl\llVert \mathbf P_{A_0}(\Delta)\bigr
\rrVert_1-\bigl\llVert \mathbf P_{A_0}^{\bot}(
\Delta)\bigr\rrVert_1 \bigr) \bigl(2Q(A_0)+\lambda
\bigl(\bigl\llVert \mathbf P_{A_0}(\Delta)\bigr\rrVert_1-
\bigl\llVert \mathbf P_{A_0}^{\bot}(\Delta)\bigr
\rrVert_1 \bigr) \bigr)
\nonumber\\[-8pt]\\[-8pt]
&&\quad=\lambda Q(A_0) \bigl\llVert \mathbf P_{A_0}(
\Delta)\bigr\rrVert_1-2\lambda Q(A_0) \bigl\llVert
\mathbf P_{A_0}^{\bot}(\Delta)\bigr\rrVert_1
\nonumber\\
&&\qquad{}+2\lambda^{2} \bigl\llVert \mathbf P_{A_0}(
\Delta)\bigr\rrVert^{2}_1+\lambda^{2} \bigl
\llVert \mathbf P_{A_0}^{\bot}(\Delta)\bigr\rrVert^{2}_1-3
\lambda^{2} \bigl\llVert \mathbf P_{A_0}(\Delta)\bigr
\rrVert_{1}\bigl\llVert \mathbf P_{A_0}^{\bot}(\Delta)
\bigr\rrVert_1.
\nonumber
\end{eqnarray}
Lemma \ref{lu1} implies that $\lambda^{2} \llVert \mathbf
P_{A_0}^{\bot}(\Delta)\rrVert^{2}_1-3\lambda^{2} \llVert
\mathbf P_{A_0}(\Delta)\rrVert_{1}\llVert \mathbf
P_{A_0}^{\bot}(\Delta)\rrVert_1\leq0$ and we obtain from
(\ref{un5})
%
\begin{eqnarray}
\label{un6}
&&\bigl(Q(\hat A_{\mathrm{SQ}})-Q(A_0) \bigr) \bigl(Q(
\hat A_{\mathrm{SQ}})+Q(A_0) \bigr)\nonumber\\[-8pt]\\[-8pt]
&&\quad\leq 4\lambda
Q(A_0) \bigl\llVert \mathbf P_{A_0}(\Delta)\bigr
\rrVert_1
-2\lambda Q(A_0) \bigl\llVert \mathbf
P_{A_0}^{\bot}(\Delta)\bigr\rrVert_1+2
\lambda^{2} \bigl\llVert \mathbf P_{A_0}(\Delta)\bigr
\rrVert^{2}_1.
\nonumber
\end{eqnarray}
Plugging (\ref{un6}) into (\ref{un1}), we get
\begin{eqnarray*}
\frac{1}{n} \bigl\llVert \Omega ( \Delta )\bigr\rrVert_2^{2}
&\leq& -2 \langle\Sigma,\Delta \rangle+4\lambda Q(A_0) \bigl\llVert
\mathbf P_{A_0}(\Delta)\bigr\rrVert_1
\\
&&{}-2\lambda Q(A_0) \bigl\llVert \mathbf
P_{A_0}^{\bot}(\Delta)\bigr\rrVert_1+2
\lambda^{2} \bigl\llVert \mathbf P_{A_0}(\Delta)\bigr
\rrVert^{2}_1.
\end{eqnarray*}
Then, by the duality between the nuclear and the operator norms, we obtain
\begingroup
\abovedisplayskip=7pt
\belowdisplayskip=7pt
\begin{eqnarray*}
\frac{1}{n} \bigl\llVert \Omega (\Delta )\bigr\rrVert_2^{2}
&\leq& 2 \llVert \Sigma\rrVert \bigl\llVert \mathbf P_{A_0}(\Delta)\bigr
\rrVert_1+2\llVert \Sigma\rrVert \bigl\llVert \mathbf
P_{A_0}^{\bot}(\Delta)\bigr\rrVert_1
\\
&&{}+4\lambda Q(A_0) \bigl\llVert \mathbf P_{A_0}(\Delta)\bigr
\rrVert_1-2\lambda Q(A_0) \bigl\llVert \mathbf
P_{A_0}^{\bot}(\Delta)\bigr\rrVert_1
\\
&&{}+2\lambda^{2} \bigl\llVert \mathbf P_{A_0}(
\Delta)\bigr\rrVert^{2}_1.
\end{eqnarray*}
Using $\lambda Q(A_0)\geq3 \llVert \Sigma\rrVert $ we compute
\[
\frac{1}{n} \bigl\llVert \Omega (\Delta )\bigr\rrVert_2^{2}
\leq\frac{14}{3}\lambda Q(A_0) \bigl\llVert \mathbf
P_{A_0}(\Delta)\bigr\rrVert_1+2\lambda^{2} \bigl
\llVert \mathbf P_{A_0}(\Delta)\bigr\rrVert^{2}_1,
\]
which leads to
\[
\frac{1}{n} \bigl\llVert \Omega (\Delta )\bigr\rrVert_2^{2}
\leq\frac{14}{3}\lambda Q(A_0) \sqrt{2\rank(A_0)}
\llVert \Delta\rrVert_2+4\lambda^{2}\rank(A_0)
\llVert \Delta\rrVert_2^{2}.
\]
The condition $4\mu m_1m_2\lambda^{2}\rank(A_0)\leq1/4$ implies that
%
\begin{equation}
\label{un7} \frac{1}{n} \bigl\llVert \Omega (\Delta )\bigr
\rrVert_2^{2}\leq\frac{14}{3}\lambda
Q(A_0) \sqrt{2\rank(A_0)}\llVert \Delta
\rrVert_2+\frac{\llVert \Delta\rrVert_2^{2}}{4\mu m_1m_2}.
\end{equation}

Set $a=\llVert \hat A_{\mathrm{SQ}}-A_0\rrVert_{\infty}$. By the
definition of $\hat A_{\mathrm{SQ}}$ we have that $ a\leq2\mathbf{a}$. We now
consider two cases, depending on whether the matrix $\frac{1}{a}
(\hat A_{\mathrm{SQ}}-A_0 )$ belongs or not to the set $ \mathcal{C}
(18\rank(A_0) )$.

\textit{Case} 1: Suppose first that $\llVert \hat A_{\mathrm{SQ}}-A_0\rrVert_{L_2(\Pi)}^{2} < a^{2}\sqrt{\frac{64 \log(d)}{\log
(6/5 ) n}}$, then (\ref{ass1}) implies that
%
\begin{equation}
\label{un8} \frac{\llVert \hat A_{\mathrm{SQ}}-A_0\rrVert^{2}_2}{m_1m_2}\leq 4\mathbf{a}^{2}\mu\sqrt{
\frac{64 \log(d)}{\log
(6/5 ) n}}
\end{equation}
and we get the statement of the Theorem \ref{thmu1} in this case.

\textit{Case} 2: It remains to consider the case $\llVert \hat
A_{\mathrm{SQ}}-A_0\rrVert_{L_2(\Pi)}^{2} \geq a^{2}\sqrt{\frac{64
\log(d)}{\log (6/5 ) n}}$. Lemma \ref{lu1} implies that
$\frac{1}{a} (\hat A_{\mathrm{SQ}}-A_0 )\in\mathcal{C}
(18\rank(A_0) )$ and we can apply Lemma \ref{thm1}. From Lemma
\ref{thm1}, (\ref{ass1}) and (\ref{un7}) we obtain that,
with probability at least $1-\frac{2}{d}$ one has
\begin{eqnarray*}
\frac{\llVert \Delta\rrVert_2^{2}}{2\mu m_1m_2}&\leq& \frac{14}{3}\lambda Q(A_0) \sqrt{2
\rank(A_0)}\llVert \Delta\rrVert_2+\frac{\llVert \Delta\rrVert_2^{2}}{4 \mu
m_1m_2}
\\
&&{}+792 a^{2}\mu m_1m_2
\rank(A_0) \bigl(\bE \bigl( \llVert \Sigma_R\rrVert
\bigr) \bigr)^{2}.
\end{eqnarray*}
A simple calculation yields
\begin{eqnarray*}
&&
\biggl(\frac{\llVert \Delta\rrVert_2}{2\sqrt{\mu
m_1m_2}}-\frac{14}{3}\lambda Q(A_0) \sqrt{2
\rank(A_0)\mu m_1m_2} \biggr)^{2}
\\
&&\quad\leq\biggl(\frac{14}{3}\lambda Q(A_0) \sqrt{2
\rank(A_0)\mu m_1m_2} \biggr)^{2}+792
a^{2}\mu m_1m_2\rank(A_0) \bigl(
\bE \bigl( \llVert \Sigma_R\rrVert \bigr) \bigr)^{2}
\end{eqnarray*}
\endgroup
and
%
\begin{eqnarray}
\frac{\llVert \Delta\rrVert_2}{2\sqrt{\mu m_1m_2}}&\leq& \frac{28}{3}\lambda Q(A_0) \sqrt{2
\rank(A_0)\mu m_1m_2}\nonumber\\[-8pt]\\[-8pt]
&&{}+\sqrt{792
a^{2}\mu m_1m_2\rank(A_0) \bigl(
\bE \bigl( \llVert \Sigma_R\rrVert \bigr) \bigr)^{2}}.\nonumber
\end{eqnarray}
This and $ a\leq2\mathbf{a}$ imply that, there exist numerical
constant $c'_1$ such that
\[
\frac{\Vert\hat A_{\mathrm{SQ}}-A_0\Vert_2^{2}}{ m_1m_2}\leq c'_1\mu^{2}
m_1m_2 \bigl( Q^{2}(A_0)
\lambda^{2}\rank(A_0)+\mathbf{a}^{2}
\rank(A_0) \bigl(\bE \bigl( \llVert \Sigma_R\rrVert
\bigr) \bigr)^{2} \bigr),
\]
which, together with (\ref{un8}), leads to the statement of the
Theorem \ref{thmu1}.

\section{\texorpdfstring{Proof of Lemma \protect\ref{l2}}{Proof of Lemma 13}}\label{pl2}
By the convexity of $Q^{2}(A)$ and using $\lambda\geq3\Delta$ we have
\begin{eqnarray*}
Q^{2}(\hat A)-Q^{2}(A_0)&\geq&-
\frac{2}{n}\sum_{i=1}^{n}
\bigl(Y_i- \langle X_i,A_0 \rangle \bigr)
\langle X_i,\hat A-A_0\rangle
\\
&=&-2 \langle\Sigma,\hat A-A_0 \rangle
\\
&\geq&-2\Vert \Sigma\Vert\Vert\hat A-A_0\Vert_1
\\
&\geq&-\frac{2}{3} \lambda\Vert\hat A-A_0
\Vert_1.
\end{eqnarray*}
Using the definition of $\hat A$, we compute
\begin{eqnarray*}
\lambda\llVert \hat A\rrVert_1-\lambda\llVert A_0
\rrVert_1&\leq& Q^{2}(A_0)-Q^{2}(\hat
A)
\\
&\leq& \tfrac{2}{3} \lambda\Vert\hat A-A_0\Vert_1.
\end{eqnarray*}

This and (\ref{ineq}) implies that
\[
\bigl\llVert \mathbf P_{A_0}^{\bot}(\hat A-A_0)
\bigr\rrVert_1\leq 5\bigl\llVert \mathbf P_{A_0}(\hat
A-A_0)\bigr\rrVert_1
\]
as stated.

\section{\texorpdfstring{Proof of Lemma \protect\ref{lu1}}{Proof of Lemma 15}}\label{plu1}

By the convexity of $Q(A)$, we have
\begin{eqnarray*}
Q(\hat A_{\mathrm{SQ}})-Q(A_0)&\geq&\frac{-(\sum_{i=1}^{n}
(Y_i-
\langle X_i,A_0 \rangle )\langle X_i,\hat
A_{\mathrm{SQ}}-A_0\rangle)/n}{Q(A_0)}
\\
&=&\frac{- \langle\Sigma,\hat
A_{\mathrm{SQ}}-A_0 \rangle}{Q(A_0)}
\\
&\geq&-\frac{\Vert\Sigma\Vert
}{Q(A_0)}\Vert\hat A_{\mathrm{SQ}}-A_0
\Vert_1
\\
&\geq&-\frac{1}{3} \lambda\Vert\hat A_{\mathrm{SQ}}-A_0
\Vert_1.
\end{eqnarray*}
Using the definition of $\hat A_{\mathrm{SQ}}$, we compute
\begin{eqnarray*}
\lambda\llVert \hat A_{\mathrm{SQ}}\rrVert_1-\lambda\llVert
A_0\rrVert_1&\leq& Q(A_0)-Q(\hat
A_{\mathrm{SQ}})
\\
&\leq& \tfrac{1}{3} \lambda\Vert\hat A_{\mathrm{SQ}}-A_0
\Vert_1.
\end{eqnarray*}
Then (\ref{ineq}) and the triangle inequality imply
\[
\bigl\llVert \mathbf P_{A_0}^{\bot}(\hat A-A_0)
\bigr\rrVert_1-\bigl\llVert \mathbf P_{A_0}(\hat
A-A_0)\bigr\rrVert_1\leq\tfrac
{1}{3} \bigl(
\bigl\llVert \mathbf P_{A_0}^{\bot}(\hat A-A_0)\bigr
\rrVert_1+\bigl\llVert \mathbf P_{A_0}(\hat
A-A_0)\bigr\rrVert_1 \bigr)
\]
and the statement of Lemma \ref{lu1} follows.
\end{appendix}

\section*{Acknowledgement}

I would like to thank Miao Weimin for his interesting comment.



\printhistory


\begin{thebibliography}{27}
%

\bibitem{chernozhukov-square}
\begin{barticle}[mr]
\bauthor{\bsnm{Belloni},~\bfnm{A.}\binits{A.}},
  \bauthor{\bsnm{Chernozhukov},~\bfnm{V.}\binits{V.}} \AND
  \bauthor{\bsnm{Wang},~\bfnm{L.}\binits{L.}}
(\byear{2011}).
\btitle{Square-root lasso: Pivotal recovery of sparse signals via conic
  programming}.
\bjournal{Biometrika}
\bvolume{98}
\bpages{791--806}.%
\bid{doi={10.1093/biomet/asr043}, issn={0006-3444}, mr={2860324}}%
\bptok{imsref}%
\end{barticle}%
\endbibitem

\bibitem{sara}
\begin{bbook}[mr]
\bauthor{\bsnm{B{\"u}hlmann},~\bfnm{Peter}\binits{P.}} \AND
  \bauthor{\bparticle{van~de} \bsnm{Geer},~\bfnm{Sara}\binits{S.}}
(\byear{2011}).
\btitle{Statistics for High-Dimensional Data:
Methods, Theory and Applications}.
\bseries{Springer Series in Statistics}.
\blocation{Heidelberg}: \bpublisher{Springer}.
\bid{doi={10.1007/978-3-642-20192-9}, mr={2807761}}
\bptok{imsref}%
\end{bbook}
\endbibitem

\bibitem{bunea}
\begin{barticle}[mr]
\bauthor{\bsnm{Bunea},~\bfnm{Florentina}\binits{F.}},
  \bauthor{\bsnm{She},~\bfnm{Yiyuan}\binits{Y.}} \AND
  \bauthor{\bsnm{Wegkamp},~\bfnm{Marten~H.}\binits{M.H.}}
(\byear{2011}).
\btitle{Optimal selection of reduced rank estimators of high-dimensional
  matrices}.
\bjournal{Ann. Statist.}
\bvolume{39}
\bpages{1282--1309}.
\bid{doi={10.1214/11-AOS876}, issn={0090-5364}, mr={2816355}}
\bptok{imsref}%
\end{barticle}
\endbibitem

\bibitem{candes-plan-noise}
\begin{barticle}[auto:STB|2013/01/18|13:50:43]
\bauthor{\bsnm{Cand{\`e}s},~\bfnm{E.~J.}\binits{E.J.}} \AND
  \bauthor{\bsnm{Plan},~\bfnm{Y.}\binits{Y.}}
(\byear{2009}).
\btitle{Matrix completion with noise}.
\bjournal{Proceedings of IEEE}
\bvolume{98}
\bpages{925--936}.
\bptok{imsref}%
\end{barticle}
\endbibitem

\bibitem{candes-recht-exact}
\begin{barticle}[mr]
\bauthor{\bsnm{Cand{\`e}s},~\bfnm{Emmanuel~J.}\binits{E.J.}} \AND
  \bauthor{\bsnm{Recht},~\bfnm{Benjamin}\binits{B.}}
(\byear{2009}).
\btitle{Exact matrix completion via convex optimization}.
\bjournal{Found. Comput. Math.}
\bvolume{9}
\bpages{717--772}.
\bid{doi={10.1007/s10208-009-9045-5}, issn={1615-3375}, mr={2565240}}
\bptok{imsref}%
\end{barticle}
\endbibitem

\bibitem{candes-tao-power}
\begin{barticle}[mr]
\bauthor{\bsnm{Cand{\`e}s},~\bfnm{Emmanuel~J.}\binits{E.J.}} \AND
  \bauthor{\bsnm{Tao},~\bfnm{Terence}\binits{T.}}
(\byear{2010}).
\btitle{The power of convex relaxation: Near-optimal matrix completion}.
\bjournal{IEEE Trans. Inform. Theory}
\bvolume{56}
\bpages{2053--2080}.
\bid{doi={10.1109/TIT.2010.2044061}, issn={0018-9448}, mr={2723472}}
\bptnote{check year}%
\bptok{imsref}%
\end{barticle}
\endbibitem

\bibitem{serbro-learning}
\begin{bmisc}[auto:STB|2013/01/18|13:50:43]
\bauthor{\bsnm{Foygel},~\bfnm{R.}\binits{R.}},
  \bauthor{\bsnm{Salakhutdinov},~\bfnm{R.}\binits{R.}},
  \bauthor{\bsnm{Shamir},~\bfnm{O.}\binits{O.}} \AND
  \bauthor{\bsnm{Srebro},~\bfnm{N.}\binits{N.}}
(\byear{2011}).
\bhowpublished{Learning with the weighted trace-norm under arbitrary sampling
  distributions. In \textit{Advances in Neural Information Processing Systems}
  (\textit{NIPS}) \textbf{24}}.
\bptok{imsref}%
\end{bmisc}
\endbibitem

\bibitem{foygel-serebro-concentration}
\begin{bmisc}[auto:STB|2013/01/18|13:50:43]
\bauthor{\bsnm{Foygel},~\bfnm{R.}\binits{R.}} \AND
  \bauthor{\bsnm{Srebro},~\bfnm{N.}\binits{N.}}
(\byear{2011}).
\bhowpublished{Concentration-based guarantees for low-rank matrix
  reconstruction. In \textit{24nd Annual Conference on Learning Theory}
  (\textit{COLT})}.
\bptok{imsref}%
\end{bmisc}
\endbibitem

\bibitem{2010arxiv10084886G}
\begin{barticle}[mr]
\bauthor{\bsnm{Ga{\"{\i}}ffas},~\bfnm{St{\'e}phane}\binits{S.}} \AND
  \bauthor{\bsnm{Lecu{\'e}},~\bfnm{Guillaume}\binits{G.}}
(\byear{2011}).
\btitle{Sharp oracle inequalities for high-dimensional matrix prediction}.
\bjournal{IEEE Trans. Inform. Theory}
\bvolume{57}
\bpages{6942--6957}.
\bid{doi={10.1109/TIT.2011.2136318}, issn={0018-9448}, mr={2882272}}
\bptnote{check year}%
\bptok{imsref}%
\end{barticle}
\endbibitem

\bibitem{gross-recovery}
\begin{barticle}[mr]
\bauthor{\bsnm{Gross},~\bfnm{David}\binits{D.}}
(\byear{2011}).
\btitle{Recovering low-rank matrices from few coefficients in any basis}.
\bjournal{IEEE Trans. Inform. Theory}
\bvolume{57}
\bpages{1548--1566}.
\bid{doi={10.1109/TIT.2011.2104999}, issn={0018-9448}, mr={2815834}}
\bptok{imsref}%
\end{barticle}
\endbibitem

\bibitem{keshavan-montanari-matrix}
\begin{barticle}[mr]
\bauthor{\bsnm{Keshavan},~\bfnm{Raghunandan~H.}\binits{R.H.}},
  \bauthor{\bsnm{Montanari},~\bfnm{Andrea}\binits{A.}} \AND
  \bauthor{\bsnm{Oh},~\bfnm{Sewoong}\binits{S.}}
(\byear{2010}).
\btitle{Matrix completion from noisy entries}.
\bjournal{J. Mach. Learn. Res.}
\bvolume{11}
\bpages{2057--2078}.
\bid{issn={1532-4435}, mr={2678022}}
\bptok{imsref}%
\end{barticle}
\endbibitem

\bibitem{keshavan-few}
\begin{barticle}[mr]
\bauthor{\bsnm{Keshavan},~\bfnm{Raghunandan~H.}\binits{R.H.}},
  \bauthor{\bsnm{Montanari},~\bfnm{Andrea}\binits{A.}} \AND
  \bauthor{\bsnm{Oh},~\bfnm{Sewoong}\binits{S.}}
(\byear{2010}).
\btitle{Matrix completion from a few entries}.
\bjournal{IEEE Trans. Inform. Theory}
\bvolume{56}
\bpages{2980--2998}.
\bid{doi={10.1109/TIT.2010.2046205}, issn={0018-9448}, mr={2683452}}
\bptok{imsref}%
\end{barticle}
\endbibitem

\bibitem{klopp-variance}
\begin{bmisc}[auto:STB|2013/01/18|13:50:43]
\bauthor{\bsnm{Klopp},~\bfnm{O.}\binits{O.}}
(\byear{2011}).
\bhowpublished{Matrix completion with unknown variance of the noise. Available
  at \url{http://arxiv.org/abs/1112.3055}}.
\bptok{imsref}%
\end{bmisc}
\endbibitem

\bibitem{klopp-rank}
\begin{barticle}[mr]
\bauthor{\bsnm{Klopp},~\bfnm{Olga}\binits{O.}}
(\byear{2011}).
\btitle{Rank penalized estimators for high-dimensional matrices}.
\bjournal{Electron. J. Stat.}
\bvolume{5}
\bpages{1161--1183}.
\bid{doi={10.1214/11-EJS637}, issn={1935-7524}, mr={2842903}}
\bptok{imsref}%
\end{barticle}
\endbibitem

\bibitem{koltchinskii-remark}
\begin{bmisc}[auto:STB|2013/01/18|13:50:43]
\bauthor{\bsnm{Koltchinskii},~\bfnm{V.}\binits{V.}}
(\byear{2011}).
\bhowpublished{A remark on low rank matrix recovery and noncommutative
  Bernstein type inequalities. In \textit{IMS Collections}, \textit{Festschrift
  in Honor of J. Wellner.}}
\bptok{imsref}%
\end{bmisc}
\endbibitem

\bibitem{koltchiskii-stflour}
\begin{bbook}[mr]
\bauthor{\bsnm{Koltchinskii},~\bfnm{Vladimir}\binits{V.}}
(\byear{2011}).
\btitle{Oracle Inequalities in Empirical Risk Minimization and Sparse Recovery
  Problems}.
\bseries{Lecture Notes in Math.}
\bvolume{2033}.
\blocation{Heidelberg}: \bpublisher{Springer}.
\bid{doi={10.1007/978-3-642-22147-7}, mr={2829871}}
\bptok{imsref}%
\end{bbook}
\endbibitem

\bibitem{koltchinskii-von}
\begin{barticle}[auto:STB|2013/01/18|13:50:43]
\bauthor{\bsnm{Koltchinskii},~\bfnm{V.}\binits{V.}}
(\byear{2011}).
\btitle{Von Neumann entropy penalization and low rank matrix estimation}.
\bjournal{Ann. Statist.}
\bvolume{39}
\bpages{2936--2973}.
\bptok{imsref}%
\end{barticle}
\endbibitem

\bibitem{Koltchinskii-Tsybakov}
\begin{barticle}[mr]
\bauthor{\bsnm{Koltchinskii},~\bfnm{Vladimir}\binits{V.}},
  \bauthor{\bsnm{Lounici},~\bfnm{Karim}\binits{K.}} \AND
  \bauthor{\bsnm{Tsybakov},~\bfnm{Alexandre~B.}\binits{A.B.}}
(\byear{2011}).
\btitle{Nuclear-norm penalization and optimal rates for noisy low-rank matrix
  completion}.
\bjournal{Ann. Statist.}
\bvolume{39}
\bpages{2302--2329}.
\bid{doi={10.1214/11-AOS894}, issn={0090-5364}, mr={2906869}}
\bptok{imsref}%
\end{barticle}
\endbibitem

\bibitem{lounici-spectral}
\begin{bmisc}[auto:STB|2013/01/18|13:50:43]
\bauthor{\bsnm{Lounici},~\bfnm{K.}\binits{K.}}
(\byear{2011}).
\bhowpublished{Optimal spectral norm rates for noisy low-rank matrix
  completion. Available at \url{http://arxiv.org/abs/1110.5346}}.
\bptok{imsref}%
\end{bmisc}
\endbibitem

\bibitem{wainwright-estimation}
\begin{barticle}[mr]
\bauthor{\bsnm{Negahban},~\bfnm{Sahand}\binits{S.}} \AND
  \bauthor{\bsnm{Wainwright},~\bfnm{Martin~J.}\binits{M.J.}}
(\byear{2011}).
\btitle{Estimation of (near) low-rank matrices with noise and high-dimensional
  scaling}.
\bjournal{Ann. Statist.}
\bvolume{39}
\bpages{1069--1097}.
\bid{doi={10.1214/10-AOS850}, issn={0090-5364}, mr={2816348}}
\bptok{imsref}%
\end{barticle}
\endbibitem

\bibitem{wainwright-weighted}
\begin{barticle}[mr]
\bauthor{\bsnm{Negahban},~\bfnm{Sahand}\binits{S.}} \AND
  \bauthor{\bsnm{Wainwright},~\bfnm{Martin~J.}\binits{M.J.}}
(\byear{2012}).
\btitle{Restricted strong convexity and weighted matrix completion: Optimal
  bounds with noise}.
\bjournal{J. Mach. Learn. Res.}
\bvolume{13}
\bpages{1665--1697}.
\bid{issn={1532-4435}, mr={2930649}}
\bptok{imsref}%
\end{barticle}
\endbibitem

\bibitem{recht-simpler}
\begin{barticle}[mr]
\bauthor{\bsnm{Recht},~\bfnm{Benjamin}\binits{B.}}
(\byear{2011}).
\btitle{A simpler approach to matrix completion}.
\bjournal{J. Mach. Learn. Res.}
\bvolume{12}
\bpages{3413--3430}.
\bid{issn={1532-4435}, mr={2877360}}
\bptok{imsref}%
\end{barticle}
\endbibitem

\bibitem{rhode-tsybakov-estimation}
\begin{barticle}[mr]
\bauthor{\bsnm{Rohde},~\bfnm{Angelika}\binits{A.}} \AND
  \bauthor{\bsnm{Tsybakov},~\bfnm{Alexandre~B.}\binits{A.B.}}
(\byear{2011}).
\btitle{Estimation of high-dimensional low-rank matrices}.
\bjournal{Ann. Statist.}
\bvolume{39}
\bpages{887--930}.
\bid{doi={10.1214/10-AOS860}, issn={0090-5364}, mr={2816342}}
\bptok{imsref}%
\end{barticle}
\endbibitem

\bibitem{serbro-collaborative}
\begin{bmisc}[auto:STB|2013/01/18|13:50:43]
\bauthor{\bsnm{Salakhutdinov},~\bfnm{R.}\binits{R.}} \AND
  \bauthor{\bsnm{Srebro},~\bfnm{N.}\binits{N.}}
(\byear{2010}).
\bhowpublished{Collaborative filtering in a non-uniform world: Learning with
  the weighted trace norm. In \textit{Advances in Neural Information Processing
  Systems} (\textit{NIPS}) \textbf{23}}.
\bptok{imsref}%
\end{bmisc}
\endbibitem

\bibitem{tropp-user}
\begin{barticle}[mr]
\bauthor{\bsnm{Tropp},~\bfnm{Joel~A.}\binits{J.A.}}
(\byear{2012}).
\btitle{User-friendly tail bounds for sums of random matrices}.
\bjournal{Found. Comput. Math.}
\bvolume{12}
\bpages{389--434}.
\bid{doi={10.1007/s10208-011-9099-z}, issn={1615-3375}, mr={2946459}}
\bptnote{check year}%
\bptok{imsref}%
\end{barticle}
\endbibitem

\end{thebibliography}
\end{document}